\documentclass[a4ja,10pt]{article}
\newtheorem{thm}{Theorem}[section]
\newtheorem{prop}[thm]{Proposition}
\newtheorem{cor}[thm]{Corollary}

\newtheorem{pf}{Proof}

\newtheorem{exmp}{Example}[section]
\newtheorem{defn}[thm]{Definition}

\usepackage{amssymb}
\usepackage{latexsym}

\newcommand{\cpn}{{{\mathbb C}{\mathbb P}^{n}}}
\newcommand{\sts}{{\mathbb C}{\mathbb P}^{3}}

\newcommand{\proj}{{\rm Proj}(S)} 

\newcommand{\oproj}{{\mathcal O}_{{\rm Proj}(S)}} 
 
\newcommand{\dbar}{d\!\!\!{\lower-0.6ex\hbox{$-$}}\!}
\newcommand{\dslash}{d\!\!\!{\lower-0.6ex\hbox{$-$}}}

\newcommand{\h}{\hbar}
\newcommand{\ott}{\lower-0.4ex\hbox{${\scriptscriptstyle{\otimes}}$}}
\newcommand{\btt}{\lower-0.2ex\hbox{${\scriptscriptstyle{\bullet}}$}}
\newcommand{\ctt}{\lower-0.2ex\hbox{${\scriptscriptstyle{\circ}}$}}
\newcommand{\dtt}{\lower-0.2ex\hbox{${\scriptscriptstyle{\diamond}}$}}
\newcommand{\odt}{\lower-0.4ex\hbox{${\scriptscriptstyle{\odot}}$}}

\begin{document}

\begin{center}
{\large\bf Quantization of Holomorphic Poisson structure \\ 
--related to Generalized K\"{a}hler structure--}
\end{center}

\par\bigskip
\begin{center}
Naoya MIYAZAKI\\
Department of Mathematics, 
Keio University, \\
Yokohama, 223-8521, JAPAN
\end{center}

\par\medskip
\noindent
{\small {\bf Abstract}:
It is known that holomorphic Poisson structures are 
closely related to theories of 
generalized K\"{a}hler geometry and bi-Hermitian structures. 
In this article, 
we introduce quantization of holomorphic Poisson structures which 
are closely related to generalized K\"{a}hler structures
/bi-Hermitian structures. 
By resulting noncommutative product $\star$ obtained via quantization, 
we also demonstrate computations with respect to concrete examples. 
\par\medskip
\noindent
{\bf Acknowledgements}: 
The present article is dedicated to 
Professor Yoshiaki Maeda on his retirement from active work 
for Keio University, 
and to Professor Hideki Omori on the occasion of his 77th birthday 
({\em Kiju}-celebration in Japan). 
The author would like to thank 
to Professors A. Asada, K. Fujii, Y. Homma, T. Suzuki, T. Tate, K. Uchino, 
and A. Yoshioka 
for their fruitful discussions and helpful comments.  

This research is partially supported 
by JSPS Grant-in-Aid for Scientific Reserch. 
}
\par\medskip\noindent
\noindent{\bf Mathematics Subject Classification (2010):} Primary 58B32; 
Secondary 53C28, 53D18, 53C55, 53D55
\par\medskip\noindent
\noindent{\bf Keywords:} Courant algebroid, 
generalized K\"{a}hler structure, deformation theory, 
Dirac structure, quantization, 
etc. 
\par\medskip\noindent

\section{Introduction} 
There are results\footnote{We recall the statement about this. 
See Theorems~2.12, 2.13 and 2.14} 
obtained in \cite{agg, gualtieri1, goto1, goto2}, 
in which we can find relations among 
holomorphic Poisson structures, 
theories of generalized K\"{a}hler geometry and bi-Hermitian structures.  
Here we first note that 
generalized K\"{a}hler geometry 
has been developed in the frame work of Courant algebroids.  
In 1990, the Courant bracket $[\![~,~]\!]$ on $TM\oplus T^*M$, 
where $M$ is a manifold, 
was introduced by T. Courant \cite{t-courant}.  
A similar bracket was defined on the double of any Lie algebroid by 
Liu, Weinstein and Xu \cite{lwx}. 
They proved that the double is not a Lie algebroid, 
but a more complicated object which is called a Courant algebroid to-day. 

The original definition of Courant algebroid in the article \cite{lwx}, 
(see also Definition~\ref{t-courant}. of the present article), 
has five conditions which are including strange anomalies 
with respect to Jacobi identity and derivation rules, etc. 
In their article, they proposed a different definition 
with a slightly modification of the original bracket
of the Courant algebroid. 
In fact, in Ph. D. thesis ``Courant algebroids, 
derived brackets and even symplectic super manifolds" 
by D. Roytenberg \cite{roytenberg}, 
he showed the equivalence of original definition and new definition. 
Relating to this topics, 
note that K. Uchino \cite{uchino} showed that two of the conditions 
and the relation (assuming the Leibniz rule) 
which are included in the original definition of Courant algebroid follow 
from the rest of the conditions. 
Moreover, descriptions of the Courant bracket via 
derived bracket was given 
(for more details, see articles Y. Kosmann-Schwarzbach \cite{k-s1} 
and D. Roytenberg \cite{roytenberg}. See also \cite{k-s2, k-sm}).

On the other hand, ``generalized complex structure" was introduced by 
N. Hitchin \cite{hitchin}. Roughly speaking, 
generalized almost complex structure ${\mathcal J}$ is defined as a section 
$SO(TM\oplus T^*M)$ satisfying ${\mathcal J}^2=-1$, where $SO$ is 
the special orthogonal group with respect to the symmetrization 
of the natural pairing $\langle~,~\rangle$ between $TM$ and $T^*M$. 
It is known that a generalized almost complex structure ${\mathcal J}$ can be 
represented ${\mathcal J}={\mathcal J}_\Phi$ 
via kernel and complex conjugate-kernel 
of a nondegenerate pure spinor $\Phi$. 
Furthermore, it is also shown that 
generalized almost complex structure ${\mathcal J}_\Phi$ 
is integrable, that is $ker~\Phi$ is 
involutive with respect to the Courant bracket, if and only if 
there exists a section $E\in \Gamma(TM\oplus T^*M)^{\mathbb C}$ such that 
$d\Phi= E\cdot\Phi$, where $\cdot$ stands for natural spinor representation 
of $CL(TM\oplus T^*M)$ on $\land^\bullet T^*M$ 
via the interior product and exterior product. 
Using frame-work Courant algebroids, 
theory of generalized complex/K\"{a}hler structure 
was developed by M. Gualtieri \cite{gualtieri1}. 
A {\em generalized K\"{a}hler structure} 
is a pare $({\mathcal J}_1, {\mathcal J}_2)$ 
of generalized complex 
(i.e. almost complex which is involutive with respect to 
the Courant bracket $[\![~,~]\!]$) structure 
s.t. ${\mathcal J}_1 {\mathcal J}_2={\mathcal J}_2 {\mathcal J}_1$ and 
$\langle-{\mathcal J}_1{\mathcal J}_2~,~\rangle$ 
is a positive definite symmetric form. 
As mentioned below, 
it is known that a generalized K\"{a}hler structure corresponds 
to a bi-Hermitian structure \cite{gualtieri1}. 
Furthermore, a complex manifold admitting bi-Hermitian structure has 
a holomorphic Poisson structure. 
Under this situation, in this article, 
we study quantization for holomorphic Poisson structures 
which appeared in generalized K\"{a}hler structures/bi-Hermitian structures. 
More precisely, we introduce a notion 
of quantization/formal symbol calculus, 
and study fundamental properties. 
By resulting noncommutative product $\star$ from quantization using 
the above Poisson structure, 
we demonstrate computations with respect to concrete examples 
for $\star$-exponential functions.  
First, using holomorphic Poisson structures,  
we establish a notion of quantization/formal symbol calculus
\footnote{
The terminology ``symbol calculus" is used in Fourier analysis 
and pseudo-differential operator 
in order to study partial differential equations
\cite{hormander-fio, hormander-I, hormander-III, shubin}. 
Especially, symbol calculus of pseudo-differential operator 
with respect to elliptic operator gives fruitful contribution 
to the index theorem. 
The notion of pseudo-differential 
operators was used to extend the class of possible deformation 
of an elliptic operator which has essential topological datas 
of the base manifold. } 
of projective schemes. 
This means that under a suitabale condition, 
we construct deformations of ring structures on projective schemes.   
We also study a concrete example, that is, 
formal symbol calculus on a projective spaces $\cpn$, 
we demonstrate concrete computations of $\star$-exponential functions. 
Precisely speaking, 
when germs of all admissible functions of holomorphic Poisson structures 
contains all germs of structure sheaf, 
we construct deformations of ring structures 
for structure sheaves on a projective scheme. 
Using holomorphic Poisson structures, as will be seen later, 
which is a typical example of Dirac structure,
we consider symbol calculus on a projective space. 
Since the structure sheaves on algebraic varieties have 
important and essential feature which plays 
crucial role to analyze their fundamental properties with respect to 
the base varieties,  
construction of symbol calculus using a holomorphic Poisson structure 
on the base variety is very interesting. 
We here state the main theorems of the present article. 

\begin{thm}\label{thm1}
Assume that $Z=[z_0:z_1:\ldots:z_n]$ is 
the homogeneous coordinate system of $\cpn$, 
and $\Lambda=\sum_{\alpha,\beta}\stackrel{\leftarrow}{\partial_{Z_{{\alpha}}}}
\Lambda^{{\alpha},{\beta}} 
        \stackrel{\rightarrow}{\partial_{Z_{{\beta}}}}$
defines a holomorphic skew-symmetric biderivation
\footnote{We use Einstein's convention unless confusing.}  
of order zero acting on the structure sheaf ${\mathcal O}_\cpn$ 
satisfying the Jacobi rule
and an assumption below: 
{ 
\begin{eqnarray}\label{lambda-relation}
\nonumber
&&
\bigr(\stackrel{\leftarrow}{\partial_{Z_{{\alpha}_1}}}
\Lambda^{{\alpha}_1,{\beta}_1} 
        \stackrel{\rightarrow}{\partial_{Z_{{\beta}_1}}}\bigr)
\cdots
\bigr(\stackrel{\leftarrow}{\partial_{Z_{{\alpha}_k}}}
\Lambda^{{\alpha}_k,{\beta}_k} 
        \stackrel{\rightarrow}{\partial_{Z_{{\beta}_k}}}\bigr)\\
&=&\stackrel{\longleftarrow}{\partial_{Z_{{\alpha}_1\ldots{\alpha}_k}}}
\Lambda^{{\alpha}_1, {\beta}_1}\cdots
\Lambda^{{\alpha}_k,{\beta}_k} 
        \stackrel{\longrightarrow}{\partial_{Z_{{\beta}_1\ldots{\beta}_k}}}, 
\end{eqnarray}
}
where 
$$
f(Z)\stackrel{\longleftarrow}{\partial_{Z_{{\alpha}_1\ldots{\alpha}_k}}}
\quad
(\mbox{resp}. ~ 
\stackrel{\longrightarrow}{\partial_{Z_{{\beta}_1\ldots{\beta}_k}}} 
g(Z))
$$  
means 
$$\partial_{Z_{\alpha_1}}\partial_{Z_{\alpha_2}}
 \cdots \partial_{Z_{\alpha_k}}f(Z)
\quad(\mbox{resp.} ~ 
\partial_{Z_{\beta_1}}\partial_{Z_{\beta_2}}
\cdots \partial_{Z_{\beta_k}}g(Z)).$$ 
Then, for any point ``$p$" and germs $[f(Z)]$ 
$~[g(Z)]$ of the stalk ${\mathcal O}_{\cpn,p}[[\natural]]$
\footnote{Here $[[\natural]]$ denotes  
either $[\mu,\mu^{-1}]]$ or $[\mu,\mu^{-1}]$, 
and we have to choose carefully in context. }, 
\begin{equation}\label{sharp-product}
\begin{array}{ll}
& [f(Z)\star g(Z)] \\
:=&[\sum_{k=0}^\infty \frac{1}{k!}\left(\frac{\mu}{2}\right)^k
\Lambda^{\alpha_1\beta_1}
\Lambda^{\alpha_2\beta_2}
\cdots \Lambda^{\alpha_k\beta_k}\\
&\qquad
\partial_{Z_{\alpha_1}}\partial_{Z_{\alpha_2}}\cdots \partial_{Z_{\alpha_k}}f(Z)\partial_{Z_{\beta_1}}\partial_{Z_{\beta_2}}\cdots \partial_{Z_{\beta_k}}g(Z)]
\end{array}
\end{equation}
defines a non-commutative and associative ring structure, 
where $\mu$ is a formal parameter. 
\end{thm}
We also have 
\begin{thm}\label{thm1-5}
Under the same assumptions and notations of Theorem~\ref{thm1}, 
the product $\star$ induces globally defined non-commutative, associative 
product on the sheaf-cohomology space  
$\sum_{k=0}^\infty H^0(\cpn,{\mathcal O}_\cpn(k))[[\natural]]$ 
where $\mu$ can be specialized a scalar (for example $\mu=1$). 
\end{thm}
Using the product $\star$, we also have 
\begin{thm}\label{thm2} 
Suppose that the same assumptions of Theorem~\ref{thm1}. 
Let $A[Z]:=ZA{}^tZ$ be a quadratic form with homogeneous degree 2. 
Then we have 
\begin{equation}\label{star-exponential}
\begin{array}{lll}
e_{\star}^{\frac{1}{\mu}A[Z]}&=&\det{}^{-1/2}
\left(\frac{e^{\Lambda A}+e^{-\Lambda A}}{2}\right)\cdot 
e^{1/\mu\left(\frac{\Lambda^{-1}}{\sqrt{-1}}\tan (\sqrt{-1}\Lambda A)\right)[Z]}\\
&\in& \sum_{k=0}^\infty H^0(\cpn,{\mathcal O}_\cpn(k)[\mu,\mu^{-1}]]).
\end{array}
\end{equation}
\end{thm}
Note that our argument is able to apply to 
weighted projective spaces, for which one can find details elsewhere. 

\section{Courant algebroid}\label{lwx}
In this section, we review fundamentals of Courant algebroid 
(See \cite{lwx}, \cite{roytenberg}). 
\subsection{Courant algebroid, Lie bialgebroid, 
and Maurer-Cartan type equation}
We start this subsection with: 
\begin{defn}\label{t-courant}
A {\em Courant algebroid} 
is a vector bundle $E\to M$ equipped with a non-degenerate 
symmetric bilinear form $(~,~)$ on the bundle, 
a skew-symmetric bracket $[\![~,~]\!]$
on the space $\Gamma(E)$ of all sections of $E$ 
and a bundle map $\rho:E\to TP$ such that 
the following properties are satisfied: 
\begin{enumerate}
\item $[\![ [\![e_1,e_2]\!],e_3]\!]+c.p.= {\mathcal D}T(e_1,e_2,e_3),~~
\forall e_1,e_2,e_3\in \Gamma(E)$,
\item $\rho([\![e_1,e_2]\!])=[\rho(e_1),\rho(e_2)], 
~~\forall e_1,e_2\in \Gamma(E)$,
\item $[\![e_1,fe_2]\!]=f[\![e_1,e_2]\!]+(\rho(e_1)f)e_2-(e_1,e_2){\mathcal D}f, 
~~\forall e_1,e_2\in \Gamma(E),\forall f\in C^\infty(M)$,
\item $\rho\circ{\mathcal D}=0$, i.e. $({\mathcal D}f,{\mathcal D}g)=0,~~
\forall f, g \in C^\infty(M)$,
\item $\rho(e)(h_1,h_2)=([\![e,h_1]\!]+{\mathcal D}(e,h_1),h_2)+
(h_1,[\![e,h_2]\!]+{\mathcal D}(e,h_2)),~~\forall e,h_1,h_2\in \Gamma(E)$,
\end{enumerate}
where ``c.p." denotes cyclic permutation, 
$T(e_1,e_2,e_3)$ is the function on the base manifold $M$ defined by 
\begin{equation}\label{trilinear}
T(e_1,e_2,e_3)=\frac{1}{3}([\![e_1,e_2]\!],e_3)+c.p., 
\end{equation}
and ${\mathcal D}:C^\infty(M)\to \Gamma(E)$ is the map defined by 
${\mathcal D}=\frac{1}{2}\beta^{-1}\rho^* d_0$, where $\beta$ is the isomorphism between $E$ and $E^*$ given bythe biliniear form. In other words, 
\begin{equation}\label{D-rho}
({\mathcal D}f,e)=\frac{1}{2}\rho(e)f . 
\end{equation} 
\end{defn}
Note that Courant algebroid gives a non-trivial example of $L_\infty$-algebra. 
(cf. \cite{roytenberg}, \cite{k}.)

Next we recall Dirac structure. 
\begin{defn}\label{Dirac}
Let $E$ be a Courant algebroid. A subbundle ${\mathcal L}$ of $E$ is called 
{\em isotoropic} if it is isotropic under the symmetric bilinear form $(~,~)$. 
${\mathcal L}$ is said to be {\em almost Dirac} if ${\mathcal L}$ 
is maximally isotropic. 
Furthermore, ${\mathcal L}$ is said 
to be {\em Dirac} if ${\mathcal L}$ is almost Dirac and 
integrable (involutive), that is, 
$\Gamma({\mathcal L})$ is closed under the bracket $[\![~,~]\!]$. 
\end{defn}
It is obvious that 
an integrable isotropic subbundle ${\mathcal L}$ of a Courant algebroid 
$(E,\rho,[\![~,~]\!],(~,~))$ is 
a Lie algebroid with $\rho|_{\mathcal L}$ and $[\![~,~]\!]$.

It is known that a Lie bialgebroid is a pair $(A,A^*)$ 
of vector bundles in duality, 
each of which is a Lie algebroid, 
such that the differential defined by 
one of them on the exterior algebra of 
its dual is a derivation of the Schouten bracket. 
Related to this definition, it is also known that 
\begin{thm}[\cite{mx}]
If a pair of Lie algebroids 
$(A,A^*)$ in duality, 
then the differential $d$ is a derivation of $(\Gamma(A^*),[~,~]_*)$ 
if and only if 
the differential $d_*$ is a derivation of $(\Gamma(A),[~,~])$. 
\end{thm}
Assume that both $A$ and $A^*$ are Lie algebroids over the base manifold $M$, 
with anchors $a$ and $a_*$ respectively. Let $E:=A\oplus A^*$. We can naturally define non-degenerate bilinear forms in the following way: 
\begin{equation}\label{bilinear}
(X_1+\xi_1, X_2+\xi_2)_{\pm}
=\frac{1}{2}(\langle \xi_1,X_2\rangle
  \pm \langle\xi_2,X_1\rangle).
\end{equation}
On $\Gamma (E)$, we can introduce a bracket by 
\begin{equation}\label{bracket}
\begin{array}{lll}
[\![e_1,e_2]\!]&=&([X_1,X_2]+L_{\xi_1}X_2-L_{\xi_2}X_1-d_*(e_1,e_2)_-) \\ 
   && +([\xi_1,\xi_2]+L_{X_1}\xi_2-L_{X_2}\xi_1+d(e_1,e_2)_-),
\end{array}
\end{equation}
where $e_1=X_1+\xi_1, e_2=X_2+\xi_2$. 

Furthermore, let $\rho:E\to TM$ be the bundle map defined by 
\begin{equation}\label{anchor}
\rho(X+\xi)=a(X)+a_*(\xi),~~\forall X \in \Gamma(A), \xi \in \Gamma(A^*).
\end{equation}
Finally, we define the operator $\mathcal D$ by 
\begin{equation}\label{derivation}
{\mathcal D}=d_*+d
\end{equation}
where $d_*:C^\infty(M)\to \Gamma(A)$ and $d:C^\infty(M)\to \Gamma(A^*)$ 
are the usual differential operator associated to Lie algebroids
\footnote{When $(A,A^*)=(\frak{g},\frak{g}^*)$, i.e. Lie bialgebra, 
the bracket reduces to the Lie bracket of Manin 
on the double $\frak{g}\oplus\frak{g}$. 
If $A=TM,~A^*=T^*M$ the the above bracket takes the form: 
$[\![X_1+\xi_1,X_2+\xi_2]\!]=[X_1,X_2]+\{L_{X_1}\xi_2-L_{X_2}\xi_1+d(e_1,e_2)_-\}$. This bracket is just known as the original Courant bracket.}.  
Under the above notations, it is known that 
\begin{thm}[\cite{lwx}]\label{Lie-bialgebroid}
If$(A,A^*)$ is a Lie bialgebroid, then $E=A\oplus A^*$  together with 
$([\![~,~]\!],\rho,(~,~)_+)$ is a Courant algebroid. 
Conversely, in a Courant algebroid 
$(E,[\![~,~]\!],\rho,(~,~))$, suppose that $L_1,L_2$ are Dirac subbundle 
transversal to each other, 
i.e. $E=L_1\oplus L_2$. Then, $(L_1,L_2)$ is a Lie bialgebroid, 
where $L_2$ is considered as the dual bundle of $L_1$ under 
the pairing $2(~,~)$.
Hence, if $(A,A^*)$ is a Lie bialgebroid, then 
$(A^*,A)$ is also a Lie bialgebroid. 
\end{thm}

Next, we recall the fundamentals of Hamiltonian operators. 
Assume that $(A,A^*)$ is a Lie algebroid. 
Suppose that $H:A^* \to A$ is a bundle map, and denote by $A_H$ 
the graph of $H$, considered as a subbundle of $E=A\oplus A^*$, 
that is, $A_H=\{H\xi+\xi|\xi\in A^*\}$. 
Under these notations, one can prove the following important 
fact: 
\begin{thm}[\cite{lwx}]\label{maurer-cartan-type}
$A_H$ is a Dirac subbundle 
(i.e. maximally isotropic, and involutive 
with respect the Courant bracket $[\![~,~]\!]$) 
if and only if 
$H$ is skew-symmetric and satisfies the following 
Maurer-Cartan type equation: 
\begin{equation}\label{maurer-cartan-type-equation}
d_*H+\frac{1}{2}[\![H,H]\!]=0, 
\end{equation}
where $H$ is considered as a section of $\land^2A$. 
\end{thm}
See Theorem 6.1, Definition 6.2 and 
Corollary 6.3 in \cite{lwx} for details. 
\par\medskip\noindent
{\bf Remark} Note that the equation (\ref{maurer-cartan-type-equation}) in 
Theorem~\ref{maurer-cartan-type} plays essential roles in order to 
discuss deformation theory of generalized K\"{a}hler structure. 
\par\bigskip\noindent
Now we recall new definitoin of Courant algebroid 
(cf. \cite{lwx} and \cite{roytenberg}). 
\begin{defn}\label{courant}
A {\em Courant algebroid} 
is a vector bundle $E\to M$ equipped with a non-degenerate 
symmetric bilinear form $(~,~)$ on the bundle, 
a (not necessarily skew-symmetric) bracket $[\![~,~]\!]^N$
on the space $\Gamma(E)$ of all sections of $E$ 
and a bundle map $\rho:E\to TP$ such that 
the following properties are satisfied: 
\begin{enumerate}
\item $[\![e_1,[\![^N e_2e_3]\!]^N]\!]^N=[\![ [\![e_1,e_2]\!]^N,e_3]\!]^N
[\![e_2,[\![e_1,e_3]\!]^N]\!]^N,~~
\forall e_1,e_2,e_3\in \Gamma(E)$,
\item $\rho([\![e_1,e_2]\!]^N)=[\rho(e_1),\rho(e_2)], 
~~\forall e_1,e_2\in \Gamma(E)$,
\item $[\![e_1,fe_2]\!]^N=f[\![e_1,e_2]\!]^N+(\rho(e_1)f)e_2, 
~~\forall e_1,e_2\in \Gamma(E),\forall f\in C^\infty(M)$,
\item $[\![e,e]\!]^N=\frac{1}{2}{\mathcal D}\langle e, e \rangle =0,~~
\forall f, g \in C^\infty(M)$,
\item $\rho(e)(h_1,h_2)=([\![e,h_1]\!]^N,h_2)+(h_1,[\![e,h_2]\!]^N),
~~\forall e,h_1,h_2\in \Gamma(E)$,
\end{enumerate}
where 
\end{defn}
\begin{prop}
Relation between original Courant bracket $[\![~,~]\!]$ and 
new Courant bracket $[\![~,~]\!]^N$ is given in the following way: 
$$[\![e_1,e_2]\!]^N=[\![e_1,e_2]\!]
+\frac{1}{2}{\mathcal D}\langle e_1, e_2 \rangle.$$
Conversely, 
$$[\![e_1,e_2]\!]=\frac{1}{2}([\![e_1,e_2]\!]^N-[\![e_2,e_1]\!]^N).$$
\end{prop}
Notice that the notion of a Dirac subbundle remains unchainged 
when switch to the new definition of a Courant algebroid. 
(See \cite{roytenberg}.) 
\par\medskip\noindent
{\bf Example} There are so many examples of Lie bialgebroids, 
and then Courant algebroids. 
For instance, Poisson manifolds, Nijenhuis manifolds, 
Poisson-Nijenhuis manifolds, 
objects constructed by solutions of classical Yang-Baxter equations, etc. 
See \cite{k-s1, k-s2, k-sm, lwx, mx}. 
\subsection{Dirac structure and generalized complex structure}
Let $M$ be a complex manifold. 
Define sheaves on $M$ by 
$$\begin{array}{l}
{T}M_p=\mbox{sheaf of germs of holomorphic vector fields around }p 
\\
{T}^*M_p=\mbox{sheaf of germs of holomorphic one forms around }p 
\end{array}$$
We use symmetric 
and skew-symmetric bilinear operations on a sheaf 
${\mathbf T}M:={T}M \oplus {T}^*M$ as 
$$
\langle (X,\xi),(Y,\eta)\rangle_{\pm}:=\frac{1}{2}\{\xi(Y)\pm\eta(X)\} 
$$
and 
$$
[\![ (X,\xi),(Y,\eta)]\!]:=
([X,Y],{L}_X\eta-L_Y\xi-i_Yd\langle(X,\xi),(Y,\eta)\rangle_-)
$$
for all $ (X,\xi),(Y,\eta)\in {\mathbf T}M.$ 
Here ${L}_X$ means the Lie derivative by $X$ and 
$i_Y$ stands for the interior product by $Y$. 
As mentioned in the previous section, 
a subbundle $D\subset {\mathbf T}M$ is called a Dirac 
structure if the following conditions are satisfied: 
\begin{enumerate}
\item[](D1)~$\langle\cdot,\cdot\rangle_+ |_D=0$;
\item[](D2)~${\rm rank}(D)$ is equal to ${\rm dim}_{\mathbf C}(M)$;
\item[](D3)~$[\![ {\mathbf T}M, {\mathbf T}M]\!]\subset {\mathbf T}M$.
\end{enumerate}
A complex manifold $M$ together with holomorphic Dirac structure 
$D\subset{\mathbf T}M$ is called a holomorphic Dirac manifold and denoted by 
$(M,D)$. 

The followings are well-known. 
\begin{exmp}
Assume that $M$ be a complex manifold 
with a holomorphic presymplectic form $\omega$. 
Then the 2-form $\omega$ induces the bundle map  
$$\omega^{\flat}:\Gamma(M)\to \Omega^1(M);~X\mapsto i_X\omega,$$
where $\Gamma(X)$ denotes the space consisting of all germs of 
holomorphic vector fields, 
and $\Omega^1(M)$ denotes the space consisting of all germs of 
holomorphic 1-differential forms. 
Then one can obtain the subbundle $\mbox{graph}(\omega^\flat)$ in 
${\mathbf T}M$ as 
$$
\mbox{\rm graph}(\omega^{\flat})_m:=
\{(X_m,i_{X_m}\omega_m)\in {T}_mM\oplus {T}^*_mM|X_m\in {T}_mM\}~~(m\in M)
$$
and can verify that $\mbox{graph}(\omega^{\flat})$ satisfies the conditions 
{\rm (D1)-(D3)}. Therefore, $(M,\mbox{graph}(\omega^{\flat}))$ defines 
a Dirac structure.
\end{exmp}
The following is also a typical example. 
\begin{exmp}
Assume that $M$ be a complex manifold 
with a holomorphic Poisson structure $\varphi$. 
Then the Poisson structure $\varphi$ induces the bundle map  
$$\varphi^{\#}:\Omega^1(M)\to\Gamma(M) 
;~\alpha \mapsto \{\bullet \to \varphi(\bullet,\alpha)\},$$ 
where $\Gamma(X)$ 
and $\Omega^1(M)$ are as above.  
Then one can obtain the subbundle $\mbox{\rm graph}(\varphi^{\#})$ in 
${\mathbf T}M$ as 
$$
\mbox{\rm graph}(\varphi^{\#})_m:=
\{(\varphi^{\#}(\i_m),\xi_m)\in {T}_mM\oplus {T}^*_mM|\xi_m\in {T}^*_mM\}~~(m\in M)
$$
and can verify that $\mbox{graph}(\varphi^{\#})$ satisfies the conditions 
{\rm (D1)-(D3)}. Therefore, $(M,\mbox{\rm graph}(\varphi^{\#}))$ defines 
a Dirac structure.
\end{exmp}
\par\bigskip
For each point $m$, a holomorphic Dirac structure $D\subset {\mathbf T}M$ 
defines two natural projections as:
$$
\rho_m:={\rm pr}_1|_{D_m}:D_m\to T_mM ~~\mbox{ and }~~
\rho_m^*:={\rm pr}_2|_{D_m}:D_m\to T_m^*M.
$$
\begin{prop}
As for the above projections, we obtain the followings: 
$${\rm ker}~\rho=D\cap (\{0\}\oplus {T}^*M)~~\mbox{ and }~~
{\rm ker}~\rho^*=D\cap ({T}M\oplus\{0\})\cdots(*1)$$
and 
$${\rm Im}~\rho= (D\cap (\{0\}\oplus {T}^*M))^\circ~~\mbox{ and }~~
{\rm Im}~\rho^*= (D\cap ({T}M\oplus\{0\}))^\circ\cdots(*2)$$
where the symbol ${}^\circ$ stands for the annihilator. 
\end{prop}

To the end of this section, we recall the fundamental facts related 
to generalized complex geometry:
Based on the original Courant bracket of the direct sum 
${\mathbf T}M=TM\oplus T^*M$ over a manifold $M$. 
The fiber bundle of the direct sum ${\mathbf T}M$ admits an 
indefinite metric $\langle~,~\rangle_+$ by which we obtain the fiber bundle 
$SO({\mathbf T}M)$ wht fiber the special orthogonal group. 
As mentioned above, an almost generalized complex structure ${\mathcal J}$ is defined as a section of the fiber bundle $SO({\mathbf T}M)$ with 
${\mathcal J}^2=-1$, which gives rise to 
the decomposition ${\mathbf T}M\otimes{\mathbb C}=
{\mathcal L}_{\mathcal J}\oplus \bar{{\mathcal L}_{\mathcal J}}$, 
where 
${\mathcal L}_{\mathcal J}$ (resp. $\bar{{\mathcal L}_{\mathcal J}}$) 
$-\sqrt{-1}$-eigenspace (resp. $\sqrt{-1}$-eigenspace).  
Almost generalized complex structures form an orbit of the action of the real Clifford group of the real Clifford algebra bundle $CL$ 
with respect to $( {\mathbf T}M,~\langle~,~\rangle_+) $. 
\begin{defn}\label{gcs}
A {\rm generalized complex structure} 
is an almost generalized complex structure 
which is integrable/involutive with respect to the original Courant bracket. 
A {\rm generalized K\"{a}hler structure} is 
a pair $({\mathcal J}_0,{\mathcal J}_1)$ consisting of commuting generalized complex structures such that $G:=-{\mathcal J}_0{\mathcal J}_1$ is 
a generalized metric. 
\end{defn}
We here recall the definition of bi-Hermitian structure: 
\begin{defn}\cite{agg}
A quadruplet $(h,J^+,J^-,b)$ is said to be a bi-Hermitian structure 
if $h$ is a Riemannian metric, $J^+,J^-$ are complex structures, and 
$b$ is a real 2-form satisfying: 
\begin{enumerate}
\item The Riemannian metric $h$ is Hermitian metric 
with respect to $J^+ (resp.J^-)$. 
\item Let $\omega_{\pm}$ be the fundamental 2-forms with respect to $J^{\pm}$. 
Then, $$d_+^c\omega_+=-d_-^c\omega_-=db,$$ 
where $\partial_{\pm},~\bar{\partial}_{\pm}$ 
are operators defined by complex structures $J^{\pm}$, 
and we put $d_{\pm}^c:=\sqrt{-1}(\bar{\partial}_{\pm}-\partial_{\pm}).$
\end{enumerate}
\end{defn}
Then we have 
\begin{thm}\cite{gualtieri1}
A bi-Hermitian structure gives a generalized K\"{a}hler structure. 
Conversely, a generalized K\"{a}hler structure gives a bi-Hermitian structure. 
\end{thm}

Furthermore, it is known that 
\begin{thm}\cite{agg}
A complex manifold admitting bi-Hermitian structure has a holomorphic 
Poisson structure. 
\end{thm}

Note that on a complex surfaces, Poisson structures are non-trivial 
section of anti-canonical bundle $K^{-1}$. 
Therefore it is well-known about holomorphic 
Poisson structures and complex Poisson surfaces. 
Moreover, it is known that 
\begin{thm}\cite{goto1,goto2}
A Poisson-K\"{a}hler manifold has non-trivial bi-Hermitian structures. 
\end{thm}
Summing up, through generalized K\"{a}hler structures, 
there are close relations between holomorphic structures and 
bi-Hermitian structures.

\section{Main results:} 
In this section, we show the main results and its super version. 
\subsection{Proofs of Theorems \ref{thm1}, \ref{thm1-5} and \ref{thm2}}
In order to give of proofs of main results, 
we begin with a slight modification of the standard theory of scheme 
Let $S=\oplus_{n=0}^\infty S_n$ be a graded commutative ring. 
Then, $S_0$ is obviously commutative and $S$ is an $S_0$-algebra. 
It is well-known that the homogeneous ideal $S_+:=\oplus_{n=1}^\infty S_n$  
is called the {\it irrelevant ideal}. Then we have  
\begin{prop}\label{noether}
A graded commutative ring $S$ is noetherian if and only if 
$S_0$ is noetherian and $S$ is finitely generated by $S_1$ as an $S_0$-algebra. \end{prop}
It is known that a projective shceme 
$$\proj :=\{\frak{p}:\mbox{ a homogeneous prime ideal}~|~ 
\neg(S_+\subset \frak{p} ) \} $$
admits the canonical scheme structure in the following way: 
Set 
$$D_+(f):=\bigr\{\frak{p} \in \proj ~|~ \neg(f\in \frak{p} ) \bigr\}, $$ 
for any homogeneous element $f\in S_d$ with degree $d$, 
then the family $\{D_+(f)\}_{f\in S_d,~d\in {\mathbb Z}_{\geq 0 }}$ forms 
a basis of open sets. Hence it gives the canonical topology $\frak{O}_{\proj}$ 
(that is, the Zariski topology) for $\proj$. 
Note that 
$\neg(f\in \frak{p} )$ means $f(\mbox{point}_{\frak{p}})\not= 0$, intuitively. 
We also set 
\begin{equation}\label{5.7}
\begin{array}{ll}
\Gamma(D_+(f),{\oproj}):=\bigr\{g/f^m~|~g\in S_m,~m\geq 0\bigr\},\\
\oproj:\frak{O}_{\proj} \ni D_+(f)  
\mapsto \Gamma(D_+(f),\oproj)\in {\rm\mathbf Mod}. 
\end{array}
\end{equation}
The functor above is the {\it structure sheaf}. 
We remark that when $g\in S_0$, we easily see that $g/1=fg/f~(fg\in S_d)$. 
Hence we may consider $g/f^m~(m\geq 1 )$ instead of $g/f^m~(m\geq 0 )$. 
We obtain that R.H.S. of (\ref{5.7}) is a part of degree $0$ of localization 
$S_f$ of $S$ by a product closed set $\{f^{\ell}\}_{\ell=0,1,2,\ldots } $. 
We denote it by $(S_f)_0$ or $S_{(f)}$. Strictly speaking,  
for any homogeneous element $f$ with ${\rm deg}(f)=d$, 
$$(S_f)_0:=S_{(f)}:=\bigr\{g/f^m~|~g\in S_{md},~m\geq 0 \bigr\}.$$
Hence summing up what mentioned above, we have 
\begin{prop}\label{scheme}
As for $D_+(f)$, 
\begin{equation}\label{affine-scheme}
(D_+(f),\oproj|_{D_+(f)})\cong {\rm Spec}(S_{(f)}).
\end{equation} 
Thus, $\proj$ is obtained by glueing of affine schemes. 
It indicates that $(\proj,\oproj)$ is a scheme in the genuin sense. 
\end{prop}
Next we consider cohomology of quasi-coherent sheaf over $\proj$. 
Assume that a graded ring $S$ is generated by $S_1$ as an $S_0$-algebra. 
For instance 
$$S=R[z_0,z_1,\ldots,z_n],~S_0=R,~S_1=\rm\{a\in S~|~{\rm deg}(a)=1\rm\}.$$
As for a quasi-coherent sheaf 
\footnote{A sheaf ${\mathcal F}$ is quasi-coherent if and only if there is a
pre-sheaf exact sequence ${\mathcal O}_U^{\oplus I}\to {\mathcal O}_U^{\oplus J}\to {\mathcal F}\to 0 $.
A sheaf ${\mathcal F}$ is coherent if and only if there is a
pre-sheaf exact sequence ${\mathcal O}_U^{\oplus n}\to {\mathcal F}\to 0~~
(n\in {\mathbb N}) $.}
${\mathcal F}$, we set 
$${\mathcal F}(m)[[\natural]]:={\mathcal F}\otimes_{\oproj} \oproj(m)[[\natural]],$$ 
and define 
\begin{equation}\label{5.21}
\Gamma_*({\mathcal F}):=
\oplus_{m\in {\mathbb Z}} \Gamma(X,{\mathcal F}(m))[[\natural]],~~
{\rm deg}(a):=m,~~(\forall a \in \Gamma(X,{\mathcal F}(m)). 
\end{equation}
Then we see that $\Gamma_*({\mathcal F}[[\natural]])$ 
is a graded $\Gamma(\oproj [[\natural]])$-module. 
For any element $f\in S_d$, we set $\alpha_d(f):=a/1$. 
Then it is well-known that 
\begin{prop}\label{thm-alpha}
The map $\alpha_d$ obtained above defines a homomorphism 
\begin{equation}\label{alpha}
\alpha_d(f):
S_d[[\natural]]\ni a \mapsto a/1 \in S(d)_{(f)}=\Gamma(D_+(f),\oproj (d) [[\natural]]).
\end{equation}
A family $\{\alpha_d(f)\}_{f\mbox{:homogeneous}}$ induces 
a module homomorphism 
\begin{equation}\label{5.22-1}
\alpha_d~:S_d[[\natural]]\to\Gamma(\proj,\oproj(d) [[\natural]]). 
\end{equation}
Hence, using the module homomorphisms $\{\alpha_d\}$, 
a graded ring homomorphism 
\begin{equation}\label{5.22}
\alpha:=\oplus_{n=0}^\infty 
\alpha_d~:S=\oplus_{n=0}^\infty S_d[[\natural]]\to\Gamma(\proj,\oproj [[\natural]])  
\end{equation}
can be defined for any quasi-coherent sheaf ${\mathcal F}$. 
Thus, $\Gamma(\proj,\oproj [[\natural]])$ admits 
a graded $S[[\natural]]$-module structure. 
\end{prop}
\begin{defn}
We denote the pair $({\rm Proj}(S),{\mathcal O}_{{\rm Proj}(S)}[[\natural]])$ 
by ${\rm Proj}(S)[[\natural]]$. 
\end{defn}
We are also interested in $\Gamma_*({\mathcal F}[[\natural]])_{(f)}$. As for any element $f\in S_d$, and $x\in\Gamma(\proj,{\mathcal F}[[\natural]])$, we see 
$x/f^n\in \Gamma_*(\proj,{\mathcal F}[[\natural]])_{(f)}.$  
We denote the restriction of $x$ to $D_+(f)$ by $x|_{D_+(f)}$. 
Then matching the degrees of $x|_{D_+(f)}$ and $(a_d(f)|_{D_+(f)})^n$ 
we get $x|_{D_+(f)}/\bigr(a_d(f)|_{D_+(f)})\bigr)^n.$ 
\begin{equation}\label{5.23-1}
\beta_{(f)}:\Gamma_*(\proj,{\mathcal F}[[\natural]])_{(f)}\ni \frac{x}{f^m}
\mapsto 
\frac{x|_{D_+(f)}}{(\alpha(f)|_{D_+(f)})^m}
\in \Gamma(D_+(f),{\mathcal F}[[\natural]]). 
\end{equation}
As similarly for $\alpha$, for any homogeneous element $g\in S_e$, 
we obtain a diagram:
\par\bigskip
$\qquad\qquad$
\unitlength 0.1in
\begin{picture}( 29.9500, 19.8500)(  4.0000,-21.9500)
\put(4.0000,-4.0000){\makebox(0,0)[lb]{$\!\!\!\!\!\!\!\!\!\!\!\!\!\!\!
\Gamma_*({\mathcal F}[[\natural]])_{(f)}$}}%
%
\special{pn 8}%
\special{pa 1000 350}%
\special{pa 3000 350}%
\special{fp}%
\special{sh 1}%
\special{pa 3000 350}%
\special{pa 2934 330}%
\special{pa 2948 350}%
\special{pa 2934 370}%
\special{pa 3000 350}%
\special{fp}%
\put(32.0000,-3.8000){\makebox(0,0)[lb]{$\Gamma(D_+(f), {\mathcal F}[[\natural]])$}}%
%
\special{pn 8}%
\special{pa 600 590}%
\special{pa 600 1990}%
\special{fp}%
\special{sh 1}%
\special{pa 600 1990}%
\special{pa 620 1924}%
\special{pa 600 1938}%
\special{pa 580 1924}%
\special{pa 600 1990}%
\special{fp}%
\put(4.0000,-22.4000){\makebox(0,0)[lb]{$\!\!\!\!\!\!\!\!\!\!\!\!\!\!\!\Gamma_*({\mathcal F}[[\natural]])_{(fg)}$}}%
%
\special{pn 8}%
\special{pa 1000 2190}%
\special{pa 3000 2190}%
\special{fp}%
\special{sh 1}%
\special{pa 3000 2190}%
\special{pa 2934 2170}%
\special{pa 2948 2190}%
\special{pa 2934 2210}%
\special{pa 3000 2190}%
\special{fp}%
\put(31.9000,-22.0000){\makebox(0,0)[lb]{$\Gamma(D_+(fg), {\mathcal F}[[\natural]])$}}%
%
\special{pn 8}%
\special{pa 3390 600}%
\special{pa 3390 2000}%
\special{fp}%
\special{sh 1}%
\special{pa 3390 2000}%
\special{pa 3410 1934}%
\special{pa 3390 1948}%
\special{pa 3370 1934}%
\special{pa 3390 2000}%
\special{fp}%
\end{picture}%

\par\bigskip\noindent
By a similar argment as above, we define an $\proj[[\natural]]$-module homomorphism 
\begin{equation}\label{5.23}
\beta_{\mathcal F}:\tilde{\Gamma_*({\mathcal F}[[\natural]])}\to {\mathcal F}[[\natural]].\end{equation}

\begin{prop}\label{serre}
Assume that a graded ring $S$ is generated by 
$S_1=\{f_1,f_2,\ldots,f_\ell\}~(\exists \ell \in {\mathbf Z}_{\geq 0})$ 
as $S_0$-algebra. Then we see 
\begin{enumerate}
\item[](i) If $S$ is a domain then the map $\alpha$ induced in 
(\ref{5.22}) is injective.
\item[](ii) If $(f_i)_{(i=1,2,\ldots,\ell)}$ are all prime ideals, then 
the map $\alpha$ is an isomorphism. 
\item[](iii) When $S=k[z_0,z_1,z_2,\ldots,z_n]$, then the map 
$\alpha$ is an isomorphism.
\end{enumerate}
\end{prop}

\begin{pf} It is almost obvious. 
However I also give a sketch of the proof for convenience. 
\par\noindent 
{\rm Injectivity of }$\alpha$. 
For any element $a \in S_m[[\mu]]$, 
if $\alpha_m(a)=0$, then we have 
\begin{equation}\label{i-inj-1}
\alpha(a)(D_+(f_i))=(a/f_i^m)\cdot(f_m^i/1)=0, 
\end{equation}
by the definition of the map $\alpha_m$. Then we see that 
$a/f_i^m=0~~(\forall i=1,2,\ldots,n)$. 
Hence there exists a positive integer $N>>0$ such that 
$f_i^Na=0~~(\forall i=1,2,\ldots,n).$ 
Since $S[[\mu]]$ is generated by $\{f_1,f_2,f_3,\ldots,f_n\}$ as an 
$S_0[[\mu]]$-algebra, and is a domain, we have $a=0$. This shows that the map 
$\alpha$ is an injective map. 
\par\medskip\noindent
{\rm Surjectivity of }$\alpha$. 
We also assume that $(f_i)$ is a prime ideal $(\forall i=1,2,,\ldots,n)$. 
For any element $h\in \Gamma(\proj,\oproj (m)[[\mu]])$, we can write 
\begin{equation}\label{i-sur-1}
h|_{D_+(f_i)} = \frac{b_i}{f_i^{m_i}}\cdot\frac{f_i^{m_i}}{1}
\in \Gamma(D_+(f_i),\oproj (m_i)[[\mu]])=S(m_i)_{(f_i)},~~b_i\in S_{m_i} . 
\end{equation}
We can assume that $m_i\geq m$. 
Note again that $D_+(f_i)\cap D_+(f_j)=D_+(f_if_j)$.

\end{pf}

We also have the following: 
\begin{prop}\label{thm5.22}
Assume that ${\mathcal F}[[\natural]]$ 
is a quasi-coherent sheaf $\proj[[\natural]]$-module. 
Then the homomorphism $\beta_{\mathcal F}$ induced 
in (\ref{5.22}) is an isomorphism. 
Furthermore, when $S=k[z_0,z_1,\ldots, z_n]$ via (\ref{5.22-1}), we see 
\begin{equation}
\begin{array}{ll}
&H^0({\mathbb P}_{k[z_0,z_1,\ldots,z_n]}^n[[\natural]])\\
:=&H^0({\rm Proj}(k[z_0,z_1,z_2,\ldots,z_n]),
{\mathcal O}_{{\rm Proj}(k[z_0,z_1,z_2,\ldots,z_n])}(m)[[\natural]] ) \\
=&
\left\{
\begin{array}{l}
0\quad\quad\quad\quad\quad\quad\quad (\mbox{if }m<0),  \\
k[z_0,z_1,z_2,\ldots,z_n]|_m \quad (\mbox{o.w.  }m\geq 0).
\end{array}\right.
\end{array}
\end{equation}
\end{prop}

\begin{defn}
Assume that $\Lambda$ is a holomorphic skew-biderivation satisfying Jacobi rule.Then, $\star$ is called {\rm symbol calculus} on 
${\rm Proj}(k[z_0,z_1,\ldots, z_n])$ if \\
$({\rm Proj}(k[z_0,z_1,\ldots, z_n])[[\natural]],\star)$ 
has an associative algebra sheaf structure 
such that 
$$
f(Z)\star g(Z) 
=f\bullet g+ \left(\frac{\mu}{2}\right)
\Lambda^{\alpha\beta}
\partial_{Z_{\alpha}}f(Z)
\partial_{Z_{\beta}}g(Z)+
\cdots, 
$$
where $f,g $ stands for germs of structure sheaf and $Z=[z_0,z_1,\ldots, z_n]$. \end{defn}

\par\bigskip\noindent
We are now in the position to give proofs of Theorems \ref{thm1}, 
\ref{thm1-5} and \ref{thm2}.~
Under the assumption (\ref{lambda-relation}), 
it is easy to check 
\begin{equation}\label{sharp-product-2}
\begin{array}{ll}
&f(Z)\sum_{k=0}^\infty 
\frac{1}{k!}\left(\frac{\mu}{2}\right)^k
\bigr(\stackrel{\leftarrow}
{\partial_{Z_{{\alpha}_1}}}
\Lambda^{{\alpha}_1,{\beta}_1} 
        \stackrel{\rightarrow}{\partial_{Z_{{\beta}_1}}}\bigr)
\cdots
\bigr(\stackrel{\leftarrow}{\partial_{Z_{{\alpha}_k}}}
\Lambda^{{\alpha}_k,{\beta}_k} 
        \stackrel{\rightarrow}{\partial_{Z_{{\beta}_k}}}\bigr)g(Z)\\
=& \sum_{k=0}^\infty \frac{1}{k!}\left(\frac{\mu}{2}\right)^k
\Lambda^{\alpha_1\beta_1}
\Lambda^{\alpha_2\beta_2}
\cdots \Lambda^{\alpha_k\beta_k}
\partial_{Z_{\alpha_1}}\partial_{Z_{\alpha_2}}\cdots \partial_{Z_{\alpha_k}}f(Z)\partial_{Z_{\beta_1}}\partial_{Z_{\beta_2}}\cdots \partial_{Z_{\beta_k}}g(Z). 
\end{array}
\end{equation}
Then the right hand side of (\ref{sharp-product-2}) 
coincides with the asymptotic expansion formula 
for product of the Weyl type pseudo-differential operators. 
Thus, it shows Theorem \ref{thm1}. 

As seen in the previous argument, 
as for sheaf cohomology of projective space, we obtain that   
\begin{equation}\label{sheaf-cohomology}
\sum_{k=0}^\infty H^0(\cpn, {\mathcal O}_\cpn(k))\cong 
\sum_{k=0}^\infty {\mathbb C}[Z]_k, 
\end{equation}
where ${\mathbb C}[Z]_k$ stands for  
the space of homogeneous polynomials of degree $k\in {\mathbf Z}_{\geq 0}$. 
Then a direct computation using (\ref{sharp-product-2}) shows 
that $\mu$ can be 
specialized a scalar. 

Finally we show Theorem \ref{thm2}.   
We would like to compute exponentials having the following form 
$f(Z)=g(t)e^{\frac{1}{\mu}Q[Z](t)}$ 
with respect to $\star$ for quadratic polynomials 
under a quite general setting.  
Let $Z=[z^1:\ldots:z^{n}]$, $A[Z]:=ZA{}^tZ$, where 
$A \in Sym(n,{\mathbb C})$, i.e. 
$A$ is an $n\times n$-complex symmetric matrix. In order to compute 
the exponential $F(t):=e_{\star}^{t\frac{1}{\mu}A[Z]}$ with respect to the 
Wyel type product formula, 
we treat the following evolution equation: 
\begin{equation}\label{evolution}
\partial_t F=\frac{1}{\mu}A[Z]\star F, 
\end{equation}
with an initial condition 
\begin{equation}\label{initial}
F_0=e^{\frac{1}{\mu}B[Z]},
\end{equation}
where $B\in Sym(n,{\mathbb C})$. 

As seen above, 
since the coefficients of bivectors are note constant, 
our setting might be seen rather different one from the situations 
considered in the appendix of the present article.
\footnote{See the article \cite{maillard} 
and in \cite{ommy1, ommy5} for original methods. 
Quillen's method employing the Cayley transform is very useful to compute superconnection character forms and 
supertrace of Dirac-Laplacian heat kernels.}. 
However, to compute exponentials, 
we can use similar methods with Cayley transform for 
homogeneous coordinates systematically, 
as will be seen below: 

Under the assumption $F(t)=g\cdot e^{\frac{1}{\mu}Q[Z]}$ 
($g=g(t),~Q=Q(t)$), we would like to find a solution of the equations 
(\ref{evolution}) and (\ref{initial}). 

Direct computations give  
{ 
\begin{eqnarray}\label{LHS}
\nonumber
\mbox{L.H.S. of } (\ref{evolution}) &=& 
g' e^{\frac{1}{\mu}Q[Z]} +g {\frac{1}{\mu}Q'[Z]}e^{\frac{1}{\mu}Q[Z]}, \\
\label{RHS}
\nonumber
\mbox{R.H.S. of } (\ref{evolution}) 
\nonumber
&\stackrel{(\ref{sharp-product})}{=}&
{\frac{1}{\mu}A[Z]}\cdot F+\frac{i\hbar}{2}
\Lambda^{i_1j_1}\partial_{i_1}{\frac{1}{\mu}A[Z]}\cdot \partial_{j_1}F  \\
\nonumber
&&\qquad-\frac{\hbar^2}{2\cdot 4} \Lambda^{i_1j_1}
\Lambda^{i_2j_2} \partial_{i_1i_2}{\frac{1}{\mu}A[Z]}
\partial_{j_1j_2}F\\
\end{eqnarray}}
where $A=(A_{ij}), \Lambda=(\Lambda^{ij})$ and $Q=(Q_{ij})$.  
Comparing the coefficient of $\mu^{-1}$ gives 
{ 
\begin{equation}
Q'[Z]=A[Z]-2{}^tA\Lambda Q[Z]-Q\Lambda A \Lambda Q[Z] . 
\end{equation}
}
Applying $\Lambda$ by left and setting 
{{} $q:=\Lambda Q$ and $a:=\Lambda A$}, 
we easily obtain 
{{} 
\begin{eqnarray}
\Lambda Q' &=& \Lambda A + \nonumber
\Lambda Q\Lambda A-\Lambda A\Lambda Q 
-\Lambda Q\Lambda A\Lambda Q \\
&=&(1+q)a(1-q).
\end{eqnarray}
}
As to the coefficient of $\mu^0$, we have 
{{} 
\begin{eqnarray}
\nonumber
g'&=&\frac{1}{2}\Lambda^{i_1j_1}\Lambda^{i_2j_2}A_{i_1i_2}gQ_{j_1j_2}  
\\
&=&-\frac{1}{2}tr(aq) \cdot g ,  
\end{eqnarray}
where $``{\it tr}"$ means the trace.  
Summing up, we have  
\begin{prop}
The equation (\ref{evolution}) is rewritten by 
{{} 
\begin{eqnarray}
\partial_t q &=& (1+q)a(1-q), \label{equation-1}\\ 
\partial_t g &=& -\frac{1}{2}tr(aq)\cdot g. \label{equation-2}
\end{eqnarray}
}
\end{prop}
In order to solve the equations (\ref{equation-1}) and (\ref{equation-2}), 
we  now recall the ``Cayley transform." 
\par\bigskip 
{{} 
\begin{prop}\label{cayley}
Set 
\begin{equation}\label{cayley-transf}
C(X):=\frac{1-X}{1+X}
\end{equation} 
if $\det (1+X)\not=0$ . Then 
\begin{enumerate}
\item $X\in sp_{\Lambda}(n,{\mathbb R})\Longleftrightarrow
\Lambda X\in Sym(n,{\mathbb R})$, \\ and then 
$C(X)\in Sp_{\Lambda}(n,{\mathbb R})$, 
where
\begin{eqnarray}
\nonumber&&Sp_{\Lambda}(n,{\mathbb R})
:=\{g\in GL(n,{\mathbb R})|{}^tg\Lambda g=\Lambda\},\\
\nonumber&&sp_{\Lambda}(n,{\mathbb R})
:=Lie(Sp_{\Lambda}(n,{\mathbb R})).
\end{eqnarray}
\item $C^{-1}(g)=\frac{1-g}{1+g}$, (the ``{\rm inverse Cayley transform}").
\item $e^{2\sqrt{-1}a}=c(-\sqrt{-1}\tan(a))$. 
\item $\log a = 2\sqrt{-1}\arctan (\sqrt{-1}C^{-1}(g))$. 
\item $\partial_t q =(1+q)a(1-q)\Longleftrightarrow 
\partial_t C(q)=-2aC(q).$ 
\end{enumerate}
\end{prop}
}

\begin{pf} 
We only show the assertion 5. 
\begin{eqnarray}\label{evolution-Cayley-0}
C(q)' 
&\stackrel{(\ref{cayley-transf})}{=}&\nonumber
\Bigr(\frac{1-q}{1+q}\Bigr)'\\
\nonumber
&\stackrel{(\ref{equation-1})}{=}&
(1+q)^{-1}(-q)'+(1+q)^{-1}(-q)'(1+q)^{-1}(1-q)\\
\nonumber&=&
-a(1-q)-a(1-q)(1+q)^{-1}(1-q)\\
\nonumber&=&-a\Bigr\{1+\frac{1-q}{1+q}\Bigr\}(1-q)\\
&\stackrel{(\ref{cayley-transf})}{=}&-2aC(q).
\end{eqnarray}
This completes the proof. 
\end{pf}

Solving the above equation 5 in Proposition \ref{cayley}, 
we have 
{{} $$C(q)=e^{-2at}C(b),$$} 
\noindent
where $b=\Lambda B$ and then 
$$
q=C^{-1}\bigr(e^{-2at}\cdot C(b)\bigr)
=C^{-1}\bigr(C(-\sqrt{-1}\tan(\sqrt{-1}at))\cdot C(b )\bigr).
$$
Hence, according to the inverse Cayley transform, 
we can get $Q$ in the following way. 
\begin{prop}
\begin{equation}\label{Q}
Q=-\Lambda\cdot C^{-1}\Bigr(C(-\sqrt{-1}\tan (\sqrt{-1} \Lambda A t)) 
\cdot C(\Lambda B)\Bigr). 
\end{equation} 
\end{prop}
Next we compute the amplitude coefficient part $g$. 
In order to find a solution, we consider the following.  
\begin{equation}
g'=-\frac{1}{2}tr(aq)\cdot g
\end{equation}
Thus, we have 
\begin{prop}
\begin{equation}
g=\det{}^{-\frac{1}{2}}
\Bigr(\frac{e^{at}(1+b)+e^{-at}(1-b)}{2}\Bigr). 
\end{equation} 
\end{prop}
\begin{pf}
First we replace 
{ 
\begin{equation}
g'=-\frac{1}{2}tr(aq)\cdot g
\end{equation}
}
by 
{ 
\begin{equation}\label{log-diff-eq}
(\log g)'=-\frac{1}{2}tr(aq).
\end{equation}
}
We also have  
{ 
\begin{eqnarray}
&&\nonumber
tr\Bigr\{\log \Bigr(\frac{e^{at}(1+b)+e^{-at}(1-b)}{2}
\Bigr)\Bigr\} '  
\\
&{=}&\nonumber
tr\Bigr\{a\frac{e^{at}(1+b)-e^{-at}(1-b)}{e^{at}(1+b)+e^{-at}(1-b)}\Bigr\}
\\
&=&\nonumber
tr(aq) . 
\end{eqnarray}
}
Combining this formula with (\ref{log-diff-eq}), we obtain  
{ 
\begin{eqnarray}
(\log g)'&=& 
\nonumber-\frac{1}{2}
tr\Bigr\{\log \Bigr(\frac{e^{at}(1+b)+e^{-at}(1-b)}{2}\Bigr)\Bigr\}'
\\
\nonumber
&=&-\frac{1}{2}\log\Bigr\{
\det \Bigr(\frac{e^{at}(1+b)+e^{-at}(1-b)}{2}\Bigr)\Bigr\}'.
\end{eqnarray}
}
This completes the proof. 
\end{pf}
Setting $t=1$, $a=\Lambda A$ and $b=0$, 
we get 
\begin{prop}\label{star-exponent}
{ 
\begin{eqnarray}
\label{star_exponential}
e_{\star}^{\frac{1}{\mu}A[Z]}
&=&\det{}^{-\frac{1}{2}}\Bigr(\frac{e^{\Lambda A}+e^{-\Lambda A}}{2}\Bigr)
\cdot e^{\frac{1}{\mu}
(\frac{\Lambda^{-1}}{\sqrt{-1}} \tan (\sqrt{-1}\Lambda A )) [Z]  } .
\end{eqnarray}
}
\end{prop}
As usual, using the $\check{\rm C}$ech resolution, 
we can compute the sheaf-cohomology 
$\sum_{k=0}^\infty H^0(\cpn, {\mathcal O}_\cpn[\mu,\mu^{-1}]])$. 
Combining it with Proposition 
\ref{star-exponent}, 
we see 
$$e_{\star}^{\frac{1}{\mu}A[Z]} \in 
\sum_{k=0}^\infty H^0(\cpn, {\mathcal O}_\cpn[\mu,\mu^{-1}]]). $$   
This completes the proof of Theorem~\ref{thm2}. 

\par\bigskip\noindent
Suppose that $\beta$ is a holomorphic Poisson structure 
Then it satisfies the Maurer-Cartan equation in Theorem~\ref{lwx} 
(cf. \cite{lwx}) 
and the adjoint action of $e^{\beta \tau}$ on ${\mathcal J}$ in the sense of 
\cite{goto2} induces an analytic family of deformations of 
generalized complex structures. 
We write it by 
${\mathcal J}_{\beta \tau} = 
\rm{Ad} (e^{\beta \tau}) {\mathcal J}$. 
Under this situation, the following was proved.  
\begin{thm}[\cite{goto2}]
Let $\beta$ be a holomorphic Poisson structure 
on a compact K\"{a}hler manifold X. 
Then we have a family of generalized K\"{a}hler structures 
denoted by $\{{\mathcal J}_{\beta \tau}, \psi_\tau\}$. 
\end{thm} 
Here we call a parameter $\tau$ a {\bf deformation parameter}. 
Then by Theorems~1.1, 1.2 and 1.3, we have the following. 
\begin{thm}\label{thm3}
Under the same assumptions as above with 
condition (\ref{lambda-relation}), 
we have a family of associative porduct $\star_\tau$ with deformation 
parameter $\tau$ by quantizing a family of holomorphic 
Poisson structures $\{\beta \tau\}$. 
Furtermore, a family of star exponentials 
$\{e_{\star_\tau}^{\frac{1}{\mu}A[Z]}\}$ exists. 
More precisely, 
\begin{eqnarray}
\label{star_exponential}
e_{\star_\tau}^{\frac{1}{\mu}A[Z]}
&=&\det{}^{-\frac{1}{2}}
\Bigr(\frac{e^{\tau \Lambda A}+e^{-\tau \Lambda A}}{2}\Bigr)
\cdot e^{\frac{1}{\mu}
(\frac{\tau^{-1}\Lambda^{-1}}{\sqrt{-1}} \tan 
(\sqrt{-1}\tau \Lambda A )) [Z]  } .
\end{eqnarray}
\end{thm}
We also obtain 
\begin{cor}
The differences of deformation parameter $\tau$ 
can be detected in the sheaf-cohomology 
$$\sum_{k=0}^\infty H^0(\cpn, {\mathcal O}_\cpn[\mu,\mu^{-1}]]).$$
\end{cor}
\par\medskip\noindent
Note that the method employed in this subsection can be extended 
to the case of the base variety is a weighted projective space 
(cf. \cite{weighted-projective-space}). 
\subsection{Super-twistor version} 
{
Consider the following diagram: 
\par\bigskip
\qquad\qquad\quad
\unitlength 0.1in
\begin{picture}( 30.0000, 23.1500)( 16.0000,-26.3000)
%
\special{pn 8}%
\special{pa 3200 800}%
\special{pa 1800 2400}%
\special{fp}%
\special{sh 1}%
\special{pa 1800 2400}%
\special{pa 1860 2364}%
\special{pa 1836 2360}%
\special{pa 1830 2338}%
\special{pa 1800 2400}%
\special{fp}%
%
\special{pn 8}%
\special{pa 3200 800}%
\special{pa 4600 2400}%
\special{fp}%
\special{sh 1}%
\special{pa 4600 2400}%
\special{pa 4572 2338}%
\special{pa 4566 2360}%
\special{pa 4542 2364}%
\special{pa 4600 2400}%
\special{fp}%
\put(32.0000,-4.0000)
{\makebox(0,0){$ \bigr((x^{\alpha,\dot{\alpha}}),[\pi_1:\pi_2]\bigr)
\in M:={\mathbb C}^{4}\times{\mathbb C}{\mathbb P}^1$}}%
\put(16.0000,-28.0000){\makebox(0,0)[lb]{$([z_1:\ldots:z_4])\in\sts$}}%
\put(44.0000,-28.0000){\makebox(0,0)[lb]{$(x^{\alpha,\dot{\alpha}})
\in{\mathbb C}^{4}$}}%
\put(44.0000,-14.0000){\makebox(0,0)[lb]{$\Pi_2$}}%
\put(16.0000,-14.0000){\makebox(0,0)[lb]{$\Pi_1$}}%
\end{picture}%

\par\vspace{1cm}
\noindent
where $x^{\alpha,\dot{\alpha}}$ are even variables. 
We set 
\begin{eqnarray}
&&\nonumber (x^{\alpha,\dot{\alpha}})
:=(x^{1,\dot{1}},x^{1,\dot{2}},x^{2,\dot{1}},x^{2,\dot{2}}),\\
&&\nonumber ([z_1:\ldots:z_4]) 
:=([x^{\alpha,\dot{1}}\pi_\alpha:x^{\alpha,\dot{2}}\pi_\alpha:\pi_1:\pi_2]). 
\end{eqnarray}
Here we use Einstein's convention
(we will often omit $\sum$ unless there is a danger of confusion). 
We call 
$([z_1:\ldots:z_4])$ 
the {\rm homogeneous coordinate system} of $\sts$. 
\begin{enumerate}
\item 
The relations\footnote{Here $[~,~]$ denotes the commutator bracket.} 
($\dot{\alpha}, \dot{\beta} = \dot{1}, \dot{2}$)  
\begin{equation}
[z^{\dot{\alpha}},z^{\dot{\beta}}]=
\hbar D^{\alpha\dot{\alpha},\beta\dot{\beta}}\pi_\alpha\pi_\beta, 
\end{equation}
where $z^{\dot{1}}:=z_1,~z^{\dot{2}}:=z_2$, 
give a globally defined non-commutative associative product
$\#$ on ${\mathbb C}{\mathbb P}^3$, 
where $\bigr(D^{\alpha\dot{\alpha},\beta\dot{\beta}}\bigr)$ 
is a skew symmetric matrix. 
\item 
Let $A[Z]$ be a homogeneous polynomial
of $z^{\dot{1}}=z_1=x^{\alpha,\dot{1}}\pi_{\alpha},~
z^{\dot{2}}=z_2=x^{\alpha,\dot{2}}\pi_{\alpha}$ 
with degree $2$. 
Then a sharp exponential function $e_{\#}^{\frac{1}{\mu}A[Z]}$ 
gives a ``function"
on ${\mathbb C}{\mathbb P}^{3}$. 
\end{enumerate} 
}
More precisely, 
\begin{thm}
Assume that $\Lambda:=\hat{\Lambda}$ 
and $A[Z]$ a homogeneous polynomial of 
$z^{\dot{1}}=x^{\alpha,\dot{1}}\pi_{\alpha},~
z^{\dot{2}}=x^{\alpha,\dot{2}}\pi_{\alpha}$ with degree $2$. 
Then a sharp exponential function $e_{\#}^{\frac{1}{\mu}A[Z]}$ 
gives a cohomology class of 
${\mathbb C}{\mathbb P}^3$ with coefficients in the sheaf 
$\sum_{k=0}^{\infty}{\cal O}_{{\mathbb C}{\mathbb P}^3}(k)$. 
\end{thm}

\par\bigskip

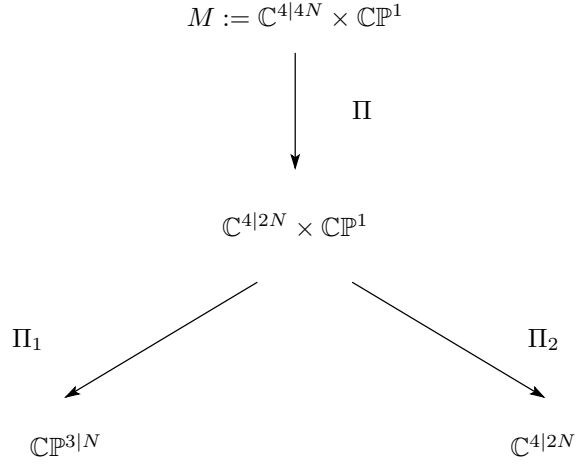
\unitlength 0.1in
\begin{picture}( 36.8500, 22.5000)(  6.1500,-31.6500)
\put(30.0000,-10.0000){\makebox(0,0){$M:={\mathbb C}^{4|4N}\times{\mathbb C}{\mathbb P}^1$}}%
%
\special{pn 8}%
\special{pa 3000 1200}%
\special{pa 3000 1800}%
\special{fp}%
\special{sh 1}%
\special{pa 3000 1800}%
\special{pa 3020 1734}%
\special{pa 3000 1748}%
\special{pa 2980 1734}%
\special{pa 3000 1800}%
\special{fp}%
\put(30.0000,-21.0000){\makebox(0,0){${\mathbb C}^{4|2N}\times{\mathbb C}{\mathbb P}^1$}}%
%
\special{pn 8}%
\special{pa 2800 2400}%
\special{pa 1800 3000}%
\special{fp}%
\special{sh 1}%
\special{pa 1800 3000}%
\special{pa 1868 2984}%
\special{pa 1846 2974}%
\special{pa 1848 2950}%
\special{pa 1800 3000}%
\special{fp}%
%
\special{pn 8}%
\special{pa 3300 2400}%
\special{pa 4300 3000}%
\special{fp}%
\special{sh 1}%
\special{pa 4300 3000}%
\special{pa 4254 2950}%
\special{pa 4254 2974}%
\special{pa 4234 2984}%
\special{pa 4300 3000}%
\special{fp}%
\put(18.0000,-32.5000){\makebox(0,0){${\mathbb C}{\mathbb P}^{3|N}$}}%
\put(43.0000,-32.5000){\makebox(0,0){${\mathbb C}^{4|2N}$}}%
\put(33.5000,-15.0000){\makebox(0,0){$\Pi$}}%
\put(16.0000,-27.0000){\makebox(0,0){$\Pi_1$}}%
\put(43.0000,-27.0000){\makebox(0,0){$\Pi_2$}}%
\end{picture}%

\par\bigskip
\noindent
where $\Pi$ denotes the chiral projection,  
we can consider non-anti-commutative deformation of super twistor space. 

In order to give a brief explanation, 
we recall the definition of super 
twistor manifold (\cite{lebr, taniguchi-miyazaki, ward-wells}). 
\begin{defn}
$(3\vert N)$-dimensional complex super manifold $Z$ is said to be 
a super twistor space if the following conditions $(1)-(3)$ are satisfied.  
\begin{enumerate}
\item[]$(1)$ $p: Z\longrightarrow {\mathbb C}{\mathbb P}^1$ 
is a holomorphic fiber bundle.
\item[]$(2)$ $Z$ has a family of holomorphic section of $p$ 
whose normal bundle is isomorphic to 
${\mathcal O}_{ {\mathbb C}{\mathbb P}^1}(1)
\oplus {\mathcal O}_{{\mathbb C}{\mathbb P}^1}(1)\oplus C^N \otimes 
\Pi {\mathcal O}_{{\mathbb C}{\mathbb P}^1}(1)$. 
\item[]$(3)$ $Z$ has an anti-holomorphic involution $\sigma$ being 
compatible with $(1),~(2)$ and $\sigma$ has no fixed point.  
\end{enumerate}
\end{defn}

We define 
${\mathbb C}{\mathbb P}^{3|N}_{*_{\alpha'}}=({\mathbb CP}^{3|N}, 
{\mathcal O}_{{\mathbb C}{\mathbb P}^{3|N}, *_\alpha'})$. 

Let $f(z\vert \xi ; \alpha^{\prime})$ be 
a local section defined in the following manner: 
\begin{equation}\label{local-represent-2}
 f(z\vert \xi ; \alpha^{\prime})=\sum_{k=0}^{N} \sum_{1\leq i_1 \leq i_2 \leq 
 \ldots \leq 
i_k \leq N} f_{i_1 i_2 \ldots i_k} (z) 
\xi^{i_1}  \xi^{i_2}  \cdots  \xi^{i_k}  
\end{equation}
where $f_{i_1\ldots i_k}(z)$ is 
a homogeneous element of $z=[z_1:z_2:z_3:z_4]$ with homogeneous degree 
$(-k)$ on ${\mathbb C}{\mathbb P}^3$. 
Then we can introduce a structure 
sheaf ${\mathcal O}_{{\mathbb C}{\mathbb P}^{3\vert N}, *\alpha^{\prime}}$
whose local section is given by $f(z\vert \xi ; \alpha^{\prime})$.  

Under these notations, we can introduce a ringed space denoted by 
${\mathbb C}{\mathbb P}^{3\vert N}_{*\alpha^{\prime}} = 
({\mathbb C}{\mathbb P}^3, 
{\mathcal O}_{{\mathbb C}{\mathbb P}^{3\vert N}, *_{\alpha'}})$. 
We shall call it non-anti-commutative complex projective super space. 

As for the non-anti-commutative deformed product $*$ 
aassociated with the non-anti-commutative complex projective super space,  
we have commutation relations of local coordinate functions: 
\begin{enumerate}
\item[]$(1)$ Let $ (z_1, z_2, \pi_1, \pi_2 \vert \xi^1,\ldots,\xi^N )$ 
be a local coordinate system of ${\mathcal P}^{3\vert N}$, 
where ${\mathcal P}^{3\vert N}$ denotes the non-anti-commutative 
open super twistor space. 
Then 
\[ \{ \xi^i ,\xi^j \}_* 
=\alpha^{\prime} C^{i\alpha, j\beta}\pi_{\alpha}\pi_{\beta}, \qquad 
(0~ \rm{o.w.})\]
\item[]$(2)$ A local coordinate system $(z_1, z_2, \pi_1, \pi_2 
\vert \xi^1,\ldots, \xi^N )$ 
of ${\mathbb C}{\mathbb P}^{3\vert N}$ satisfies 
\[ \{ \xi^i ,\xi^j \}_* 
=\alpha^{\prime} C^{i\alpha, j\beta}\pi_{\alpha}\pi_{\beta}, \qquad 
(0~\rm{o.w.}) \]
\end{enumerate}
These arguments leads us quantization of body-part 
and soul-part of super-twister spaces. 
Here we do not explain more the notion and notations which appeared above 
and do not give the proof of them. 
For details, see \cite{taniguchi-miyazaki}. 

\section{Appendix: Riccati-type equation appeared in 
$*$-transcendental elements }
\label{PBW-theorem} Star-exponential/transcendental elements 
are important objects in deformation quantization
(cf. \cite{bffls, dl, fedosov, k, s, y}).
In the present section, for convenience, 
we explain about transcendentally extended Weyl algebras 
which gave a crew for our study. 
See \cite{maillard, ommy0, ommy, ommy1, ommy2, ommy3, ommy4}, 
about relating topics. 

Let ${\mathfrak S}(n)$ (resp. ${\mathfrak A}(n)$) be 
the spaces of complex symmetric matrices 
(resp. skew-symmetric matrices). Note that 
${\mathfrak M}(n){=}
{\mathfrak S}(n)\oplus{\mathfrak A}(n)$.
For an arbitrary fixed $n{\times}n$-complex matrix 
$\Lambda{\in}{\mathfrak M}(n)$, 
we define a product ${*}_{_{\Lambda}}$ on the space of polynomials   
${\mathbb C}[{} u]$ by the formula 
\begin{equation}
 \label{eq:KK}
\begin{array}{lll}
\nonumber f*_{_{\Lambda}}g&=&fe^{\frac{i\h}{2}
(\sum\overleftarrow{\partial_{u_i}}
{\Lambda}_{ij}\overrightarrow{\partial_{u_j}})}g\\
&=&\sum_{k}\frac{(i\h)^k}{k!2^k}
{\Lambda}_{i_1j_1}\!{\cdots}{\Lambda}_{i_kj_k}
\partial_{u_{i_1}}\!{\cdots}\partial_{u_{i_k}}f\,\,
\partial_{u_{j_1}}\!{\cdots}\partial_{u_{j_k}}g. 
\end{array}   
\end{equation}
$\partial_{u_{i}}$ acts as a 
derivation of the algebra 
$({\mathbb C}[{} u],*_{_{\Lambda}})$ 
in the biderivation 
$\sum\overleftarrow{\partial_{u_i}}
{\Lambda}_{ij}\overrightarrow{\partial_{u_j}}$. 
Then $({\mathbb C}[{} u],*_{_{\Lambda}})$ is an associative algebra. 
Remark that it is not necessary commutative, in general. 
If the matrix $\Lambda=({\Lambda}_{ij})$ 
is symmetric, then the obtained algebra 
is commutative and it is isomorphic  
to the standard polynomial algebra with $\h$. 

Note that, for any other constant symmetric matrix $K$,
we can also define a new product $*_{_{\Lambda+K}}$ by the formula    
\begin{equation}\label{general-product-formula}
\begin{array}{lll}
f*_{_{\Lambda{,}K}}g&=&f
e^{\frac{i\h}{2}
(\sum\overleftarrow{\partial_{u_i}}
{K}_{ij}{*_{_{\Lambda}}}\overrightarrow{\partial_{u_j}})}g\\
&=&\sum_{k}\frac{(i\h)^k}{k!2^k}
{K}_{i_1j_1}\cdots{K}_{i_kj_k}
(\partial_{u_{i_1}}\cdots\partial_{u_{i_k}}f){*_{_{\Lambda}}}
(\partial_{u_{j_1}}\cdots\partial_{u_{j_k}}g)\\
&=&
\sum_{k}\frac{(i\h)^k}{k!2^k}
{(\Lambda{+}K)}_{i_1j_1}
\cdots{(\Lambda{+}K)}_{i_kj_k}
\partial_{u_{i_1}}\cdots\partial_{u_{i_k}}f
\partial_{u_{j_1}}\cdots\partial_{u_{j_k}}g, 
\end{array}
\end{equation}
where we used the following formula:
\begin{equation}
  \label{eq:Hochsch}
\begin{array}{lll}
&&e^{\frac{i\h}{4}\sum{K}_{ij}\partial_{u_i}\partial_{u_j}}
\Big(\big(e^{-\frac{i\h}{4}\sum{K}_{ij}
\partial_{u_i}\partial_{u_j}}f\big)
{*_{_\Lambda}}
\big(e^{-\frac{i\h}{4}\sum{K}_{ij}\partial_{u_i}\partial_{u_j}}g\big)\Big)\\
&=&fe^{\frac{i\h}{2}(\sum\overleftarrow{\partial_{u_i}}
{*_{\Lambda}}{K}_{ij}{*_{_\Lambda}}
\overrightarrow{\partial_{u_j}})}g= f{*}_{_{\Lambda{+}K}}g.
 \end{array}
\end{equation}
Before discussing the next topics, we fix notations. 
Set $\Lambda=K{+}J$ where $K$, $J$ are the symmetric part
and the skew part of $\Lambda$, respectively. 
Then the commutator is
$[u_i,u_j]={i\h}J_{ij}$. 
We also use 
\begin{equation}\label{Weyl}
{{} u}=(u_1,u_2,\cdots,u_{2m})=
({{}}{{} u},{{}}{{} v}),\quad  
{{}}{{} u}=({{}}{u}_1,\cdots,{{}}{u}_m),\,\,
{{}}{{} v}=({{}}{v}_1,\cdots,{{}}{v}_m).
\end{equation}
and use the standard skew-symmetric 
matrix 
$J=\left[
{\footnotesize 
{\begin{array}{cc} 
   0 & {-}I\\ 
   I & 0 
 \end{array}}} 
\right]$.

Note that according to the choice of $K=0, K_0, {-}K_0$ 
where 
$$
(0,\,\, K_0, {-}K_0)= 
\left(
 \left[
{\footnotesize
{\begin{array}{cc}
   0 & 0\\
   0 & 0
\end{array}}}
 \right],\,\,
\left[
{\footnotesize
{\begin{array}{cc}
   0 & I\\
   I & 0
 \end{array}}}
\right],\,\,
\left[
{\footnotesize
{\begin{array}{cc}
   0 &\!\!{-}I\\
   {-}I&\!\!0
 \end{array}}}\right]\,\,
\right). 
$$
When $K=0$ (resp. $K_0$, $K_0$), the ordering is called Weyl ordering 
(resp. normal ordering, anti-normal ordering). 
For each ordered expression, the  product formulas are given  
respectively by the following formulas:}
\begin{equation}\label{ppformula}
 \begin{array}{lll}
f({{} u}){*{_{_0}}}g({{} u})&=&
f\exp 
\frac{\h i}{2}\{\overleftarrow{\partial_{v}} 
     {\wedge}\overrightarrow{\partial_{u}}\}g,
     \quad{\mbox{(Moyal product formula)}}\\
f({{} u}){*{_{_{K_0}}}}g({{} u})&=&
f\exp {\h i}\{\overleftarrow{\partial_{v}}\,\, 
       \overrightarrow{\partial_{u}}\}g,
       \qquad{\mbox{($\Psi$DO-product formula)}}  \\
f({{} u}){*{_{_{{-}K_0}}}}g({{} u})&=&
f\exp{-\h i}\{\overleftarrow{\partial_{u}}\,\, 
       \overrightarrow{\partial_{v}}\}g,
       \quad{\mbox{($\overline{\Psi}$DO-product formula)}}   
 \end{array}
\end{equation} 
where 
$\overleftarrow{\partial_{v}}{\wedge}
\overrightarrow{\partial_{u}}
=\sum_i(\overleftarrow{\partial_{{{}}{v}_i}}
\overrightarrow{\partial_{{{}}{u}_i}}
-\overleftarrow{\partial_{{{}}{u}_i}}
\overrightarrow{\partial_{{{}}{v}_i}})$ and 
$\overleftarrow{\partial_{v}}\,\, 
       \overrightarrow{\partial_{u}}
=\sum_i\overleftarrow{\partial_{{{}}{v}_i}}\,\, 
       \overrightarrow{\partial_{{{}}{u}_i}}$.
\par\bigskip 
For $f({} u)\in H{\!o}l({\mathbb C}^{2m})$, 
the direct calculation via the product formula 
(\ref{eq:KK}) by using Taylor expansion gives the following:

\begin{equation}\label{extend}
\begin{array}{lll}
&&e^{s\frac{1}{i\h}\langle{{} a},{{} u}\rangle}
{*_{_K}}f({{} u})=
e^{s\frac{1}{i\h}\langle{{} a},{{} u}\rangle}
f({{} u}{+}\frac{s}{2}{{} a}(K{+}J)),\\ 
&&
f({{} u}){*_{_K}}
e^{-s\frac{1}{i\h}\langle{{} a},{{} u}\rangle}=
f({{} u}{+}\frac{s}{2}{{} a}(-K{+}J))
e^{-s\frac{1}{i\h}\langle{{} a},{{} u}\rangle}
\end{array} 
\end{equation} 
as natural extension of the product formula. This also gives the 
associativity of computations involving two functions 
of exponential growth and a holomorphic function. 

We use notation ${:}{\bullet}{:}_{_K}$ 
which stands for the ordering ($K$-ordered expression parameter) for 
elements of Weyl algebra $W_{\h}(2m)$. 
For instance, we write 
$$
{:}u_i{*}u_j{:}_{_K}{=}u_iu_j{+}\frac{i\h}{2}(K{+}J)_{ij}, \quad  etc.
$$
Set 
$H_*=a{{}} u_*^2{+}b{{}} v_*^2{+}2c{{}} u{\ctt}{{}} v$, 
${{}} u{\ctt}{{}} v=\frac{1}{2}({{}} u{*}{{}} v{+}{{}} v{*}{{}} u)$, and  
$c^2{-}ab=D$. It is easy to see that 
${:}H_*{:}_0=a{{}} u^2{+}b{{}} v^2{+}2c{{}} u{{}} v$.
We would like to see the form of Weyl ordered expression for 
${:}e_*^{t(a{{}} u^2+b{{}} v^2+2c{{}} u{\ctt}{{}} v)}{:}_{0}$. 
For the purpose, we set 
${:}e_*^{t(a{{}} u^2+b{{}} v^2+2c{{}} u{\ctt}{{}} v)}{:}_0= F(t,{{}} u,{{}} v)$ and consider the real analytic solution of the evolution equation 
\begin{equation}
  \label{eq:siki}
\frac{\partial}{\partial t}F(t,{{}} u,{{}} v)=
(a{{}} u^2\!+\!b{{}} v^2\!+\!2c{{}} u{{}} v){*_0}F(t,{{}} u,{{}} v), \quad 
F(0,{{}} u,{{}} v)=1.
\end{equation}
By the product formula (\ref{general-product-formula}),  
we have 
$$
\begin{array}{lll}
&&
(a{{}} u^2\!+\!b{{}} v^2\!+\!2c{{}} u{{}} v){*_0}
F(t,{{}} u,{{}} v)\\
&=&
(a{{}} u^2\!+\!b{{}} v^2\!+\!2c{{}} u{{}} v)F 
+{\hbar i}\{(b{{}} v\!+\!c{{}} u)\partial_{{{}} u}F-(a{{}} u\!+\!c{{}} v)\partial_{{{}} v}F\}\\ 
&&
{}\qquad\qquad-\frac{\hbar^2}{4}\{b\partial_{{{}} u}^2F\!-\!
2c\partial_{{{}} v}\partial_{{{}} u}F\!+\!a\partial_{{{}} v}^2F\}. 
\end{array}
$$
By using a function $f(x)$ of one variable   
$$
{:}e_*^{t(a{{}} u^2+b{{}} v^2+2c{{}} u{{}} v)}{:}_0= 
f_t(a{{}} u^2{+}b{{}} v^2{+}2c{{}} u{{}} v), 
$$ 
we get a simplified form
$$
\begin{array}{ll}
&(a{{}} u^2{+}b{{}} v^2{+}2c{{}} u{{}} v){*_0}f_t(a{{}} u^2{+}b{{}} v^2{+}2c{{}} u{{}} v)\\
=&(a{{}} u^2{+}b{{}} v^2{+}2c{{}} u{{}} v)f_t(a{{}} u^2{+}b{{}} v^2{+}2c{{}} u{{}} v)\\
   &\qquad\qquad
   -{\hbar^2}(ab{-}c^2)(f'_t(a{{}} u^2{+}b{{}} v^2{-}2c{{}} u{{}} v)\\
     &\qquad\qquad\qquad\qquad
        +f''_t(a{{}} u^2{+}b{{}} v^2{+}2c{{}} u{{}} v)
             (a{{}} u^2{+}b{{}} v^2{+}2c{{}} u{{}} v)), 
\end{array}
$$
where $x=a{{}} u^2+b{{}} v^2+2c{{}} u{{}} v$. Then we obtain the equation 
(we call this Riccati-type equation)
\begin{equation}\label{eq:diff-eq-sq}
\frac{d}{dt}f_t(x) =xf_t(x)+{\hbar^2}D(f'_t(x)+xf''_t(x))
\end{equation}
where $D=c^2-ab$ is the discriminant of $H_*$.
%
Set 
$f_t(x)=g(t)e^{h(t)x}$. Substituting this into 
(\ref{eq:diff-eq-sq}), we obtain  
$$
\big\{g'(t)-{D\hbar^2}g(t)h(t)+ 
 xg(t)\{h'(t)-1-{D\hbar^2}h(t)^2\}\big\}e^{h(t)x}=0.
$$
Hence, $h'(t)-1-{D\hbar^2}h(t)^2=0$. Thus, we obtain 
$$
h(t)=\frac{1}{\hbar\sqrt{D}}\tan(\hbar(\sqrt{D})t).
$$ 

Next, solving  
$$
g'(t)-g(t)D\hbar^2 
   \frac{1}{\hbar\sqrt{D}}\tan(\hbar(\sqrt{D})t)=0 , 
$$
we have 
$g(t)= \frac{1}{\cos(\hbar(\sqrt{D})t)}$.
Summing up what mentioned above, we have 
the solution of the differential equation 
(\ref{eq:diff-eq-sq}) with the initial function $1$ in the following way:
\begin{thm}\label{sisuu}
The Weyl ordered expression of the 
$*$-exponential function 
$e_*^{t(a{{}} u^2{+}b{{}} v^2{+}2c{{}} u{\ctt}{{}} v)}$ is given by 
$$
{:}e_*^{t(a{{}} u^2{+}b{{}} v^2{+}2c{{}} u{\ctt}{{}} v)}{:}_0=
\frac{1}{\cos(\hbar{\sqrt{D}\,t})}
 \exp\left(\frac{1}{\hbar\sqrt{D}}
\tan(\hbar{\sqrt{D}\,t})(a{{}} u^2{+}b{{}} v^2{+}2c{{}} u{{}} v)\right) 
$$
where  
$\frac{1}{\hbar\sqrt{D}}\tan(\hbar{\sqrt{D}\,t}){=}t$ 
in the case $D{=}0$. 
\end{thm}

\end{document}